\documentstyle[sprocl]{article}

\def\Journal#1#2#3#4{{#1} {\bf #2}, #3 (#4)}

\def\JPA{{\em J. Phys.} A}
\def\MPLA{{\em Mod. Phys. Lett.} A}
\def\IJMPA{{\em Int. J. Mod. Phys.} A}
\def\LMP{\em Lett. Math. Phys.}

\def\tm{\tilde m}
\def\tiota{\tilde{\iota}}
\def\tDelta{\tilde{\Delta}}
\def\tepsilon{\tilde{\epsilon}}
\def\tS{\tilde S}
\def\tone{\tilde 1}
\newcommand{\nudot}{\stackrel{\nu}{\cdot}}

\begin{document}

\title{
COLOURED HOPF ALGEBRAS AND THEIR DUALS
}
\author{C. QUESNE \footnote{Directeur de recherches FNRS}}

\address{Universit\'e Libre de Bruxelles, Brussels, Belgium}


\maketitle

\abstracts{Coloured Hopf algebras, related to the coloured Yang-Baxter equation,
are reviewed, as well as their duals. The special case of coloured quantum
universal
enveloping algebras provides a coloured extension of Drinfeld and Jimbo
formalism. The universal $\cal T$-matrix is then generalized to the coloured
context, and shown to lead to an algebraic formulation of the coloured
$RTT$-relations, previously proposed by Basu-Mallick as part of a coloured
extension of Faddeev, Reshetikhin, and Takhtajan approach to quantum groups and
quantum algebras.}
%
%
Some years ago, Kundu and Basu-Mallick~\cite{kundu,basu} extended the
Faddeev, Reshe\-tikhin, Takhtajan (FRT) $RLL$- and $RTT$-relations by replacing
the $R$-matrix, which is a solution of the Yang-Baxter equation (YBE) by a
solution
$R^{\lambda,\mu}$ of the coloured YBE, where $\lambda$, $\mu$ are non-additive
(multicomponent) parameters, referred to as `colour' indices. In such a
substitution, both $L$- and $T$-matrices acquire a colour parameter
dependence.\par
%
%
More recently, we introduced coloured Hopf algebras and coloured universal $\cal
R$-matrices, and obtained as a by-product a coloured extension of the
Drinfeld and Jimbo (DJ) approach to quantum universal enveloping algebras
\linebreak (QUEA's) of Lie algebras.\cite{cq96}\par
%
%
Here, we show that the coloured extensions of the FRT and DJ formulations of
quantum groups and quantum algebras can be related provided the notions of
dually
conjugate Hopf algebras and universal $\cal T$-matrix~\cite{fronsdal} are
generalized to the coloured context.\cite{cq97a}\par
%
%
A coloured Hopf algebra ${\cal H}^c = \left({\cal H}, {\cal C}, {\cal
G}\right) =
\bigl({\cal H}_q, +, m_q, \iota_q, \Delta^{\lambda,\mu}_{q,\nu},
\epsilon_{q,\nu},
S^{\mu}_{q,\nu};$ $k, {\cal Q}, {\cal C}, {\cal G}\bigr)$ is built from (i)
a set ${\cal
H} = \{\,{\cal H}_q \mid q\in {\cal Q}\,\}$ of standard Hopf
algebras~${\cal H}_q$
(over some field $k = {\bf R}$ or $\bf C$), whose parameters~$q$ vary in a
parameter set~$\cal Q$, (ii) a colour set~$\cal C$, whose elements are colour
parameters $\lambda$, $\mu$, $\nu$,~$\ldots$, and (iii) a colour group
${\cal G} =
\{\,\sigma^{\nu} \mid \nu \in {\cal C}\,\}$, whose elements $\sigma^{\nu}: {\cal
H}_q \to {\cal H}_{q^{\nu}}$ are algebra isomorphisms.\cite{cq96} It is endowed
with ordinary multiplication~$m_q$, and unit~$\iota_q$, but coloured
comultiplication~$\Delta^{\lambda,\mu}_{q,\nu}$, counit~$\epsilon_{q,\nu}$, and
antipode~$S^{\mu}_{q,\nu}$. For instance, $\Delta^{\lambda,\mu}_{q,\nu}:
{\cal H}_{q^{\nu}} \to {\cal H}_{q^{\lambda}} \otimes {\cal H}_{q^{\mu}}$
is defined
by $ \Delta^{\lambda,\mu}_{q,\nu} \equiv \left(\sigma^{\lambda} \otimes
\sigma^{\mu}\right) \circ \Delta_q \circ \sigma_{\nu}$ in terms of the standard
comultiplication~$\Delta_q$ of~${\cal H}_q$, and of $\sigma^{\lambda},
\sigma^{\mu}, \sigma_{\nu} \equiv \left(\sigma^{\nu}\right)^{-1} \in {\cal
G}$.\par
%
%
Whenever the starting Hopf algebras~${\cal H}_q$ are quasitriangular, the
coloured
Hopf algebra~${\cal H}^c$ is characterized by the existence of a coloured
universal
$\cal R$-matrix, denoted by~${\cal R}^c$, and defined as the set ${\cal
R}^c \equiv
\{\, {\cal R}^{\lambda,\mu}_q \mid q \in {\cal Q}, \lambda, \mu \in {\cal C}
\,\}$.\cite{cq96} Its elements ${\cal R}^{\lambda,\mu}_q \equiv
\left(\sigma^{\lambda} \otimes \sigma^{\mu}\right) \left({\cal R}_q\right) \in
{\cal H}_{q^{\lambda}} \otimes {\cal H}_{q^{\mu}}$, where ${\cal R}_q$ is the
universal $\cal R$-matrix of ${\cal H}_q$, satisfy the coloured YBE.\par
%
%
Provided the ${\cal H}_q$'s have Hopf duals~${\cal H}_q^*$, i.e., there
exist some
doubly nondegenerate bilinear forms $\langle\, , \,\rangle_q$ with respect to
which ${\cal H}_q$ and ${\cal H}_q^*$ are paired in the usual
way,\cite{fronsdal} a
coloured Hopf dual of~${\cal H}^c$ can be defined, and is denoted by ${\cal
H}^{c*} =
\left({\cal H}^*, {\cal C}, {\cal G}\right) = \left({\cal H}_q^*, +,
\tm^{\nu}_{q,\lambda,\mu}, \tiota^{\nu}_q, \tDelta_q, \tepsilon_q,
\tS^{\nu}_{q,\mu}; k, {\cal Q}, {\cal C}, {\cal G}\right)$.\cite{cq97a}
Here ${\cal
H}^* = \{\,{\cal H}_q^* \mid q\in {\cal Q}\,\}$, $\cal C$ has the same
meaning as
before, while $\cal G$ is now isomorphically represented by the set
$\{\,\rho^{\nu}
\mid \nu \in {\cal C}\,\}$, whose elements $\rho^{\nu}: {\cal H}_q^* \to {\cal
H}_{q^{\nu}}^*$, defined by $\left\langle \rho^{\nu}(x), \sigma^{\nu}(X) \right
\rangle_{q^{\nu}} = \left\langle x, X \right\rangle_q$ for any $x \in {\cal
H}_q^*$,
$X \in {\cal H}_q$, $q \in {\cal Q}$, and~$\nu \in {\cal C}$, are coalgebra
isomorphisms. The coalgebra maps $\tDelta_q$, $\tepsilon_q$ of the coloured Hopf
dual~${\cal H}^{c*}$ coincide with those of~${\cal H}_q^*$, whereas its algebra
maps $\tm^{\nu}_{q,\lambda,\mu}$, $\tiota^{\nu}_q$, and antipode
$\tS^{\nu}_{q,\mu}$ are coloured, and dual to $\Delta^{\lambda,\mu}_{q,\nu}$,
$\epsilon_{q,\nu}$, and $S^{\mu}_{q,\nu}$, respectively. For instance, the
coloured
multiplication $\tm^{\nu}_{q,\lambda,\mu}: {\cal H}_{q^{\lambda}}^* \otimes
{\cal
H}_{q^{\mu}}^* \to {\cal H}_{q^{\nu}}^*$ is defined by
$\tm^{\nu}_{q,\lambda,\mu}
\equiv \rho^{\nu} \circ \tm_q \circ \left(\rho_{\lambda} \otimes \rho_{\mu}
\right)$ in terms of the standard multiplication~$\tm_q$ of~${\cal H}_q^*$,
and of
$\rho^{\nu}, \rho_{\lambda} \equiv \left(\rho^{\lambda}\right)^{-1},
\rho_{\mu} \equiv \left(\rho^{\mu}\right)^{-1} \in {\cal G}$. An alternative notation
is $x(\lambda) \nudot y(\mu) \equiv \tm^{\nu}_{q,\lambda,\mu}(x \otimes y)$,
where $x \in {\cal H}_{q^{\lambda}}^*$, and $y \in {\cal H}_{q^{\mu}}^*$.\par
%
%
In the special case where the ${\cal H}_q$'s are QUEA's of Lie algebras
$U_q(g)$,
${\cal H}^c = U^c(g)$ is called coloured QUEA, and provides a coloured
extension of
DJ~formalism.\cite{cq96} The Hopf duals~${\cal H}_q^*$ (if they exist) are then
Hopf algebras of quantized functions on the corresponding Lie groups,
$Fun_q(G) =
G_q$, and the coloured Hopf dual is denoted by $G^c$.\cite{cq97a}\par
%
%
The universal $\cal T$-matrix of $G_q$,\cite{fronsdal} defined by ${\cal T}_q =
\sum_A x^A \otimes X_A$, where $\{x^A\}$ and $\{X_A\}$ denote dual bases
of~$G_q$ and~$U_q(g)$, respectively, can be generalized into a coloured $\cal
T$-matrix of~$G^c$, ${\cal T}^c \equiv \{\, {\cal T}^{\lambda}_q \mid q \in
{\cal Q},
\lambda \in {\cal C} \,\}$, where ${\cal T}^{\lambda}_q \equiv {\cal
T}_{q^{\lambda}}.$\cite{cq97a} The elements of the latter satisfy various
relations
extending well-known properties of~${\cal T}_q$. In particular, if $U^c(g)$
has a
coloured universal $\cal R$-matrix, then ${\cal T}^c$ of~$G^c$ fulfils the
relation
\begin{equation}
  \left(\tone_{q^{\nu}} \otimes {\cal R}^{\lambda,\mu}_q\right)
  \left({\cal T}^{\lambda}_{q,1\mu} \nudot {\cal T}^{\mu}_{q,2\lambda}\right) =
  \left({\cal T}^{\mu}_{q,2\lambda} \nudot {\cal T}^{\lambda}_{q,1\mu}\right)
  \left(\tone_{q^{\nu}} \otimes {\cal R}^{\lambda,\mu}_q\right).
\label{eq:colRTT}
\end{equation}
Here
\begin{eqnarray}
  {\cal T}^{\lambda}_{q,1\mu} & \equiv & \sum_A x^A \otimes (X_A \otimes
        1_{q^{\mu}}) \in G_{q^{\lambda}} \otimes U_{q^{\lambda}}(g) \otimes
         U_{q^{\mu}}(g), \\
  {\cal T}^{\mu}_{q,2\lambda} & \equiv & \sum_B x^B \otimes
         (1_{q^{\lambda}} \otimes X_B) \in G_{q^{\mu}} \otimes
U_{q^{\lambda}}(g)
         \otimes U_{q^{\mu}}(g),
\end{eqnarray}
where ${\cal T}^{\lambda}_q = \sum_A x^A \otimes X_A$, ${\cal T}^{\mu}_q =
\sum_B x^B \otimes X_B$ are the universal $\cal T$-matrices of
$G_{q^{\lambda}}$,
$G_{q^{\mu}}$, and $1_{q^{\lambda}}$, $1_{q^{\mu}}$, $\tone_{q^{\nu}}$ are the
units of $U_{q^{\lambda}}(g)$, $U_{q^{\mu}}(g)$, $G_{q^{\nu}}$,
respectively.\par
%
%
In the defining $n$-dimensional matrix representation $D_q$ of~$U_q(g)$, the
elements of the coloured universal $\cal R$- and $\cal T$-matrices are
represented by $n^2 \times n^2$, and $n \times n$ matrices, $R^{\lambda,\mu}_q
\equiv \left(D_{q^{\lambda}} \otimes D_{q^{\mu}}\right) \left({\cal
R}^{\lambda,\mu}_q\right)$, $T^{\lambda}_q \equiv \sum_A x^A
D_{q^{\lambda}}(X_A)$, respectively. Equation~(\ref{eq:colRTT}) then becomes the
matrix equation
\begin{equation}
  R^{\lambda,\mu}_q \left(T^{\lambda}_{q1} \nudot T^{\mu}_{q2}\right) =
  \left(T^{\mu}_{q2} \nudot T^{\lambda}_{q1}\right) R^{\lambda,\mu}_q,
\end{equation}
which is essentially identical with the coloured $RTT$-relations, introduced by
Basu-Mallick.\cite{basu}\par
%
%
It is straightforward to extend the formalism presented here to deal with graded
algebraic structures.\cite{cq97a,cq97b} One then gets coloured Hopf
superalgebras in general, and coloured QUEA's of Lie superalgebras in
particular.
In the latter case, the coloured Hopf dual is built from the Hopf
superalgebras of
quantized functions on the corresponding Lie supergroups. Quasitriangularity
implies that the elements of ${\cal R}^c$ satisfy the coloured graded YBE.\par
%
%
The general construction of coloured pairs $\left(U^c(g), G^c\right)$ has been
illustrated by giving some explicit results for the two-parameter
deformations of
$\bigl(U(gl(2)), Gl(2)\bigr)$, and $\bigl(U(gl(1/1)), Gl(1/1)\bigr)$.\cite{cq97a}\par
%
%
\section*{References}

\end{document}